\documentclass[12pt]{iopart}

\usepackage{iopams}  
\usepackage{graphics}
\newtheorem{theorem}{Theorem}[section]

\newtheorem{lemma}[theorem]{Lemma}

\newtheorem{assumption}[theorem]{Assumption}
\newtheorem{remark}[theorem]{Remark}
\renewcommand{\theequation}{\thesection.\arabic{equation}}  
\begin{document}
\title[Stability result for a time dependent potential in a waveguide]{Stability result for a time dependent potential in a waveguide}
\author{Patricia Gaitan$^1$ and Yavar Kian$^2$}
\address{$^1$ Universit\'e d'Aix-Marseille, IUT Aix-en-Provence
et LATP, UMR CNRS 7353}
\address{$^2$ Universit\'e d'Aix-Marseille, IUT Aix-en-Provence et CPT, UMR CNRS 6205}
\ead{patricia.gaitan@univ-amu.fr and yavar.kian@univ-amu.fr}
\begin{abstract}We consider the operator $H:= \partial_t -\Delta+V$ 
in $2$D or $3$D waveguide.
With an adapted global Carleman estimate with singular weight functions we give a stability result for the time dependent part of the potential for this particular geometry.
Two cases are considered: the bounded waveguide with mixed Dirichlet and Neumann conditions and the open waveguide with Dirichlet boundary conditions.
\end{abstract}

\maketitle
%
\section{Introduction}
%
\setcounter{equation}{0}
We first consider a bounded waveguide 
$\Omega =(-L,L)  \times \mathcal{D}$ in $\mathbb{R}^d$
with $d=2$ or $d=3$. In the two-dimensional case, $ \mathcal{D}=(0,h)$, where $h$ is a fixed positive constant, while in the three-dimensional case, $ \mathcal{D}$ is a connected, bounded and open domain of $\mathbb{R}^2$ with $\mathcal{C}^{\infty}$ boundary denoted $\Gamma_1$.
We denote $\Gamma_1=\Gamma_1^+ \cup \Gamma_1^-$.\
We will denote $x=(x_1,x_2)$ a generic point of $\Omega$ where $x_1 \in(-L,L)$ and $x_2 \in \mathcal{D}$.
We consider the heat equation
\begin{equation}\label{systub}
\left \{ \begin{array}{ll}
\partial_t u-\Delta u +V(t,x)u=0, &\mbox{ in }(0,T)\times\Omega,\\
u(t,x)=b(t,x), &\mbox{on } (0,T)\times [-L,L] \times \partial\mathcal{D} ,\\
\partial_{\nu} u(t,\pm L,x_2)=k_\pm(t,x_2) ,
 &\mbox{ on } (0,T)\times\mathcal{D} ,\\
u(0,x)=u_0(x), &\mbox{ on } \Omega,
\end{array}\right.
\end{equation}
where $V(t,x)=q(t,x_2)f(x_1) $ and $f(x_1) >0$.
The aim of this paper is to give a stability and uniqueness result for the time dependent part of the potential
$q(t,x_2)$ using  global Carleman estimates.
We denote by $\nu$ the outward unit normal to $\Omega$ on $\partial \Omega$. 
\\ \noindent
We shall use the following notations $Q= (0,T)\times\Omega$ and
$\Sigma=(0,T)\times\partial\Omega$.\\ \\ \noindent
Our problem can be stated as follows:\\  \noindent
Is it possible to determine the coefficient $q(t,x_2)$ from the measurement of
$\partial_{\nu}(\partial_{x_1} u)$ on $(0,T) \times(-L,L) \times \Gamma_1^+$, where $\Gamma_1^+$ is a part of $\partial \mathcal{D}$?\\ \\ \noindent
Let $u$ (resp. $\widetilde{u}$) be a solution of (\ref{systub}) associated with
($q$, $f$, $b$, $u_0$) (resp. ($\widetilde{q}$, $f$, $b$, $u_0$)).

\noindent
Our main result is
$$\|q-\widetilde{q}\|^2_{L^2((\varepsilon,T-\varepsilon)\times \mathcal{D})} \leq C_\varepsilon\|\partial_{\nu}(\partial_{x_1} u)-
\partial_{\nu}(\partial_{x_1}\widetilde{u})\|^2_{L^2((0,T)\times (-L,L) \times\Gamma_1^+)},\quad 0<\varepsilon<\frac{T}{2},$$
where $C_\varepsilon$ is a positive constant which depends on $(\Omega, \Gamma_1^+, \varepsilon, T)$
and where the above norms are weighted Sobolev norms. 
\vskip 0.5cm
 \noindent
Using a method introduce in \cite{AT}, we first derive a global Carleman estimate with singular weight
for the operator $H:=\partial_t-\Delta$ with a boundary term on 
a part $\Gamma$ of the boundary $\Gamma$ of $\mathcal{D}$.
Then using these estimate and following the method developed by
Choulli and Yamamoto \cite{CY06}, we give
a stability and uniqueness result for the time dependent part  $q(t,x_2)$ of the potential $V(t,x)$.
\\
\noindent
For the first time, the method of Carleman estimates was introduced in the field
of inverse problems in the work of Bukhgeim and Klibanov \cite{BK}. A recent book
by Klibanov and Timonov \cite{KT} is devoted to the Carleman estimates 
applied to inverse coefficient problems.
\\ \noindent
The problem of recovering time independent coefficients has attracted considerable attention recently
and  many theoretical results exist. Regarding time dependent coefficient few results exist.
 In the case of source term, Canon and Esteva \cite{CE} established uniqueness and \textit{a priori} estimates for the heat conduction equation with over specified data. Choulli and Yamamoto \cite{CY06} obtained a stability result, in a restricted class, for the inverse problem of determining a source term $f(x,t)$ from Neumann boundary data for the heat equation in a bounded domain. In a recent work, Choulli and Yamamoto \cite{CY11} considered the inverse problem of  finding a control parameter $p(t)$ 
that reach a desired temperature $h(t)$ along a curve $\gamma(t)$ for a parabolic semi-linear equation with homogeneous Neumann boundary data and they established existence, uniqueness as well as Lipschitz stability. Using optic geometric solution, Choulli \cite{Ch09} considered the inverse problem of determining a general time dependent coefficient of order zero for parabolic equations  from Dirichlet to Neumann map and he proved uniqueness as well as stability. In \cite{E07} and \cite{E08}, Eskin considered the same inverse problem for hyperbolic and the Schr\"odinger equations with time-dependent electric and magnetic potential and he established uniqueness by gauge invariance. \\
The idea introduce in \cite{Y} allows us to take into account the 
particular geometry of our domain. Indeed, the $x_1$-derivative do not alter the Dirichlet condition.\\
\noindent
In a second part, we will consider an open waveguide 
$\Omega =\mathbb{R}  \times \mathcal{D}$ in $\mathbb{R}^d$.
We will denote $x=(x_1,x_2)$ a generic point of $\Omega$ where $x_1 \in \mathbb{R}$ and $x_2 \in \mathcal{D}$.
We consider the heat equation
\begin{equation}\label{systu}
\left \{ \begin{array}{ll}
\partial_t u-\Delta u +V(t,x)u=0 \;\mbox{ in }\; (0,T)\times\Omega,\\
u(t,x)=b(t,x) \;\mbox{ on }\; (0,T)\times\partial\Omega ,\\
u(0,x)=u_0(x) \;\mbox{ on }\; \Omega,
\end{array}\right.
\end{equation}
where $V(t,x)=q(t,x_2)f(x_1) $ and $f(x_1) \ge c_{min}>0$.
Using an approach similar to the previous case we will give a stability and uniqueness result for the time dependent part of the potential
$q(t,x_2)$ using  global Carleman estimates.
We denote by $\nu$ the outward unit normal to $\Omega$ on $\partial \Omega=\mathbb{R} \times \partial\mathcal{D}$.\\ 
\noindent
This paper is organized as follows.  In section $2$, using an adapted global Carleman estimate for the 
operator $H$, we give a stability result for the coefficient $q$ of problem (\ref{systub}).
In section $3$, we consider the open wave guide and we establish a stability result for the coefficient $q$ of problem (\ref{systu}).
%
\section{Stability result for a bounded waveguide}
%
In this section we consider problem (\ref{systub}). We will establish a stability
result and deduce a uniqueness  for the
coefficient $q$. We give the result for the bounded waveguide with mixed boundary   conditions. The Carleman estimate  will be the key
ingredient in the proof of such a  stability estimate.

\noindent
 From now on we set \[\Omega=(-L,L)\times \mathcal{D},\  Q=(0,T)\times\Omega,\  \Sigma_1=(0,T)\times[-L,L]\times\partial \mathcal{D},\] \[\Sigma'_2=(0,T)\times\{-L\}\times \mathcal{D},\  \Sigma''_2=(0,T)\times\{L\}\times \mathcal{D},\  \Sigma_2=\Sigma'_2\cup\Sigma''_2\] and \[\Sigma=(0,T)\times\partial\Omega.\]
 %
\subsection{Carleman estimate}
%
Let us introduce the differential operator 
\[P=\partial_t-\Delta_x+A(t,x)\cdot\nabla_x\cdot+B(t,x)\]
with $A\in L^\infty(Q,\mathbb{R}^d)$ and $B\in L^\infty(Q)$. Let $\Gamma^+_1$ be a closed subset of $\partial \mathcal{D}$, let $\alpha\in(-L,L)$ and let $\Gamma^+$ be defined by
\[\Gamma^+=(-L,L)\times\Gamma^+_1.\]
 Let $\psi_2$ be a 
${\mathcal C}^4(\mathbb{R}^{d-1})$ function satisfying the following conditions:
\begin{enumerate}\label{psi2}
\item $\psi_2(x_2) >0$ in $\overline{\mathcal{D}}$,
\item  There exists $C_0>0$ such that $|\nabla {\psi_2}| \geq C_0>0 \;\;\mbox{ in } \;\;\mathcal D$,
\item ${\partial}_{\nu} {{\psi_2}}\leq 0\;\;\mbox{on}\;\;(\partial \mathcal D\setminus\Gamma^+_1)$.
\end{enumerate}
For the proof of existence of  a function satisfying these conditions, we refer to \cite{FI} and \cite{CIK}. Now let $\psi_1\in{\mathcal C}^4(\mathbb{R})$ be a function satisfying the following conditions
 \begin{enumerate}\label{psi1}
 \item  $\psi_1(x_1)>0$ for $x\in(-L,L)$,
\item  $\psi_1'(x_1)<0$ for $x\in(-L,\alpha)$,
\item  $\psi_1'(x_1)>0$ for $x\in(\alpha,L)$,
\item  $\psi_1'(-L)=\psi_1'(L)=0$.
\end{enumerate}

One can easily prove existence of a function satisfying these conditions. Choose $\psi(x_1,x_2)=\psi_1(x_1)\psi_2(x_2)$. Then, $\psi$ is a ${\mathcal C}^4(\mathbb{R}^{d})$  function satisfying the conditions
\begin{assumption}We have:
\label{psi}\begin{itemize} 
\item $\psi(x) >0$ in $\overline{\Omega}$,
\item  There exists $C_0>0$ such that $|\nabla {\psi}| \geq C_0>0 \;\;\mbox{ in } \;\;\Omega$,
\item ${\partial}_{\nu} {{\psi}}\leq 0\;\;\mbox{on}\;\;(\partial \Omega\setminus\Gamma^+)$,
\item  $\partial_{x_1}\psi(x)<0$ for $x\in(-L,\alpha)\times \mathcal{D}$,
\item  $\partial_{x_1}\psi(x)>0$ for $x\in(\alpha,L)\times \mathcal{D}$.
\end{itemize}\end{assumption}

Now, let us introduce the function
\begin{equation}\label{ca1}\eta(t,x)=g(t)\left(e^{2\lambda|{\psi}|_\infty}-e^{\lambda \psi(x)}\right),\quad \rho>0\end{equation}
with
\[g(t)=\frac{1}{t(T-t)}.\]
We consider the following Carleman estimate.
\begin{theorem}\label{t5} \emph{( Theorem 3.4, \cite{Ch09})} 
Let Assumption \ref{psi} be fulfilled. Then, there exist three constants $s_0$, $C$ and $\lambda$ depending of $\Omega$, $T$, $\Gamma^+$, $|{A}|_{L^\infty(Q,\mathbb{R}^d)}$ and $|{B}|_{L^\infty(Q)}$ such that
\begin{eqnarray}\label{ca2} 
\int_Qe^{-2s\eta}\left[(s g)^{-1}(\Delta u)^2+(s g)^{-1}(\partial_t u)^2+s g|{\nabla u}|^2+(s g)^3 u^2\right] dx\ dt\\ \nonumber
\leq C\left(\int_Qe^{-2s\eta}(Pu)^2\ d x\ d t+\int_{(0,T)\times\Gamma^+}e^{-2s \eta}s g(\partial_\nu u)^2 d\sigma\ dt\right)
\end{eqnarray}
for $s\geq s_0$ and $u\in \mathcal C^{2,1}(Q)$, $u=0$ on $\overline{\Sigma}$.
\end{theorem} 
%
\subsection{Inverse Problem}
%
Let $u$ be solution of
$$\left \{ \begin{array}{lll}
   \partial_t u -\Delta u+q(t,x_2)f(x_1)u=0 & \mbox{in} &Q,\\
   u(t,x)=b(t,x)  & \mbox{on} & \Sigma_1,\\
   \partial_{\nu} u(t,-L,x_2)=k_-(t,x) ,\;\;
\partial_{\nu} u(t,L,x_2)=k_+(t,x) &\mbox{on} &(0,T)\times\mathcal{D} ,\\
   u(0,x)=u_0(x) & \mbox{in} & \Omega,
\end{array}\right.$$
and $\widetilde{u}$ be solution of
\begin{equation}\label{eqi}\left \{ \begin{array}{lll}
  \partial_t \widetilde{u} -\Delta \widetilde{u}+\widetilde{q}(t,x_2)f(x_1)\widetilde{u}  =0 & \mbox{in} & Q,\\
   \widetilde{u}(t,x)=b(t,x)  & \mbox{on} &\Sigma_1,\\
   \partial_{\nu} \widetilde{u}(t,-L,x_2)=k_-(t,x_2) ,\;\;
\partial_{\nu} \widetilde{u}(t,L,x_2)=k_+(t,x_2)&\mbox{on} &(0,T)\times\mathcal{D} ,\\
   \widetilde{u}(0,x)=u_0(x) & \mbox{in} & \Omega,
\end{array}\right.\end{equation}
%
Let us consider the following conditions.

\begin{assumption} \label{as1} Here we assume that 
\begin{itemize}
\item $q(t,x_2)f(x_1), \widetilde{q}(t,x_2)f(x_1)\in\mathcal C^{1+\alpha,\frac{\alpha}{2}}(\overline{Q}),$
\item $b\in\mathcal C^{3+\alpha,1+\frac{\alpha}{2}}([0,T]\times[-L,L]\times\partial\mathcal D),$
\item $u_0\in \mathcal C^{3,\alpha}(\overline{\Omega}),$
\item $k^\pm\in\mathcal C^{2+\alpha,1+\frac{\alpha}{2}}([0,T]\times\overline{\mathcal D}),$
\item $\partial_tb(0,x)-\Delta_xu_0(x)+q(0,x_2)f(x_1)u_0(x)=0,\quad x=(x_1,x_2)\in[-L,L]\times\partial\mathcal D,$
\item $\partial_tb(0,x)-\Delta_xu_0(x)+\widetilde{q}(0,x_2)f(x_1)u_0(x)=0,\quad x=(x_1,x_2)\in[-L,L]\times\partial\mathcal D,$
\item  $f>0$, $b>0$,  $u_0>0$,
\item $\widetilde{u}(t,-L,x_2)>0,\ \ \widetilde{u}(t,L,x_2)>0,\quad (t,x_2)\in [0,T]\times \overline{\mathcal D}.$
\end{itemize}
\end{assumption}
\noindent
Notice that Assumption \ref{as1} and the maximum principle implies that $\widetilde{u}(t,x)>0$ for $(t,x)\in\overline{Q}$. Moreover, according to \cite[Chapter 4 Section 5]{LSU}, we have $u$, $\widetilde{u}$, $\partial_{x_1}u$, $\partial{x_1}\widetilde{u}\in \mathcal C^{2,1}(\overline{Q})$.\\
\noindent
\begin{remark}
 Assume that $b(t,x)$ can be extended to a function $b_1\in\mathcal C^{3+\alpha,1+\frac{\alpha}{2}}(\overline{\Sigma})$ such that $b_1>0$ and  the compatibility condition
\[\partial_tb_1(0,x)-\Delta_xu_0(x)+\widetilde{q}(0,x_2)f(x_1)u_0(x)=0,\quad x\in\Gamma\]  is fulfilled. Let $w$ be the solution of 
$$\left \{ \begin{array}{lll}
   \partial_t w -\Delta w+\widetilde{q}(t,x_2)f(x_1)w=0 & \mbox{in} &Q,\\
    w(t,x)=b_1(t,x)  &\mbox{on} & \Sigma,\\
   w(0,x)=u_0(x) & \mbox{in} & \Omega,
\end{array}\right.$$
Choose $k_\pm$ such that 
\[ \partial_{\nu} w(t,\pm L,x_2)=k_\pm(t,x_2),\quad (t,x_2)\in[0,T]\times\mathcal{D}.\]
Then $w$ will be a solution of \emph{(\ref{eqi})} and by uniqueness we obtain
\[w(t,x)=\widetilde{u}(t,x),\quad (t,x)\in Q.\]
It follows that $\widetilde{u}_{\vert \overline{\Sigma}}=b_1>0$ and the last item of Assumption \ref{as1} can
 be removed.
\end{remark}
Now, if we set $v=u- \widetilde{u}$, then $v$ satisfies 
%
$$\left \{ \begin{array}{lll}
 \partial_t v -\Delta v+q(t,x_2)f(x_1)v=(\widetilde{q}(t,x_2)-q(t,x_2))f(x_1)\widetilde{u}
& \mbox{in} & Q,\\
   v(t,x)=0  & \mbox{on} & \Sigma_1,\\
   \partial_{\nu} v(t,-L,x_2)=0 ,\;\;
\partial_{\nu} v(t,L,x_2)=0 &\mbox{on} &(0,T)\times\mathcal{D} ,\\
   v(0,x)=0 & \mbox{in} & \Omega,
\end{array}\right.$$
%
Thus with the change of function $\displaystyle{w=\frac{v}{f\widetilde{u}}}$, $w$ is solution of the following system
$$\left \{ \begin{array}{lll}
   \partial_t w -\Delta w+\mathbb{A}\cdot \nabla w +aw=\widetilde{q}(t,x_2)-q(t,x_2)
 & \mbox{in} &Q,\\
   w(t,x)=0  & \mbox{on} &\Sigma_1,\\
   \partial_{\nu} w(t,-L,x_2)=0 ,\;\;
\partial_{\nu} w(t,L,x_2)=0&\mbox{on} &(0,T)\times\mathcal{D} ,\\
   w(0,x)=0 & \mbox{in} & \Omega,
\end{array}\right.$$
%
where
$$\mathbb{A}=\frac{-2}{f\widetilde{u}}\nabla(f\widetilde{u}) \ \ \mbox{and}\ \ 
a=\frac{\partial_t(f\widetilde{u})-\Delta(f\widetilde{u})+V}{f\widetilde{u}}.$$
\noindent
We consider the $x_1$-derivative of the previous system and we set $z:=\partial_{x_1} w$.
Let us observe that for all function $g\in\mathcal C^1(\overline{Q})$ we have
\begin{equation}\label{inv1}\partial_\nu g(t,L,x_2)=\partial_{x_1}g(t,L,x_2),\quad x_2\in\mathcal D\end{equation}
and
\begin{equation}\label{inv2}\partial_\nu g(t,-L,x_2)=-\partial_{x_1}g(t,-L,x_2),\quad x_2\in\mathcal D.\end{equation}
Moreover, if $g(t,x)=0$ for $(t,x)\in(-L,L)\times\partial \mathcal D$, since $\partial_{x_1}$ is a tangent derivative on $(-L,L)\times\partial D$, we have
\begin{equation}\label{inv3}\partial_{x_1} g(t,x_1,x_2)=0,\quad (t,x_1,x_2)\in(0,T)\times(-L,L)\times\partial \mathcal D.\end{equation}
From (\ref{inv1}), (\ref{inv2}), (\ref{inv3}) and the fact that $z\in\mathcal C^{2,1}(\overline{Q})$, we deduce that 
 $z$ is the solution of
\begin{equation} \label{syst1-z}
\left \{ \begin{array}{lll}
   \partial_t z -\Delta z+\mathbb{A}\cdot \nabla z +az+B_1z=B_2\partial_{x_2}w +bw
    & \mbox{in} & Q,\\
   z(t,x)=0  & \mbox{on} & \Sigma,\\
   z(0,x)=d(x) & \mbox{in} & \Omega
\end{array}\right.
\end{equation}
with 
$$B_1:=-2\partial_{x_1}\left(\frac{\partial_{x_1}(f\widetilde{u})}{f\widetilde{u}}\right),\ \ B_2:=2\partial_{x_1}\left(\frac{\partial_{x_2}(f\widetilde{u})}{f\widetilde{u}}\right)    \ \ \mbox{and}\ \ b=-\partial_{x_1} a.$$
Note that Assumption \ref{as1} implies that
$\mathbb{A}$, $B_1$, $B_2$ and $a$, $b$ are bounded. 

\noindent
 Set 
 \begin{equation}\label{inv4} I_1(z)=\int_Qe^{-2s\eta}\left[(s g)^{-1}(\Delta z)^2+(s z)^{-1}(\partial_t z)^2+s g|{\nabla z}|^2+(s g)^3 u^2\right]\ d x\ d t.\end{equation}
 Applying the Carleman estimate (\ref{ca2}) to $u=z$,  we obtain
%
$$I_1(z) \leq C\left[
  \int_{(0,T)\times\Gamma^+}e^{-2s \eta}s g(\partial_\nu z)^2 d\sigma dt \right.$$
$$ \left.  + \int_{Q} e^{-2s \eta}\ (| \nabla z |^2+ | z |^2)\ d x \ d t+
 \int_{Q} e^{-2s \eta}\ (| \nabla w |^2+ | w |^2)\ d x \ d t\right].$$
  %
  The second integral of the right hand side of the previous estimate is "absorbed"
  by the left hand side, for $s$ sufficiently large. For the last integral, we need the following lemma proved in \cite{Ch09}:
  \begin{lemma} \label{l1}
  Let $F$ be a function in $\mathcal C(\overline Q))$. Then we have the following estimate:
$$\int_Q \left |  \int_{\alpha}^{x_1} F(t,\xi,x_2) d\xi \right | ^2 e^{-2s\eta} dx_1\; dx_2\; dt \leq C
   \int_Q \left |  F(t,x)  \right | ^2 e^{-2s\eta} dx\; dt .$$
  \end{lemma}
 { \bf{Proof}}\\ 
 We recall here the proof of this lemma. An application of the Cauchy-Schwarz inequality yields
    \begin{equation}\label{l1a} 	\hspace{-2cm}\int_Q \left |  \int_{\alpha}^{x_1} F(t,\xi,x_2) d\xi \right | ^2 e^{-2s\eta} dx_1\; dx_2\; dt
    \leq C \int_Q     \int_{\alpha}^{x_1}  | {F(t,\xi,x_2)} | ^2 d\xi    e^{-2s\eta} dx_1\; dx_2\; dt.\end{equation}
  Let $r(x_1,\xi)$ be defined by
  $r(x_1,\xi)=e^{-2s[\phi(t,x_1,x_2)-\phi(t,\xi,x_2)]}$.
  Note that
  \[\partial_\xi r(x_1,\xi)=-2s\lambda\frac{\partial\psi}{\partial x_1}(\xi,x_2) g(t)e^{\lambda \psi(\xi,x_2)}r(x_1,\xi).\]
  Thus, Assumption \ref{psi} implies that $\partial_\xi r(x_1,\xi)<0$ for $\alpha<\xi<x_1<L$ and $\partial_\xi r(x_1,\xi)>0$ for $-L<x_1<\xi<\alpha$. Then, in the region 
  \[\{(x_1,\xi):\ \alpha\leq\xi\leq x_1\leq L\}\cup \{(x_1,\xi):\ -L\leq x_1\leq\xi\leq\alpha\}\]
  we have
  \begin{equation}\label{l1b}r(x_1,\xi)\leq r(x_1,x_1)=1.\end{equation}
Applying estimate (\ref{l1a}) , we get
$$ \int_Q \left |  \int_{\alpha}^{x_1} F(t,\xi,x_2) d\xi \right | ^2 e^{-2s\eta} dx\; dt $$
$$\leq C\left[\int_0^T\int_{\mathcal D}\int_{-L}^\alpha\int^{\alpha}_{x_1}  | {F(t,\xi,x_2)} | ^2 e^{-2s\eta(t,\xi,x_2)}r(x_1,\xi)d\xi\; dx_1\; dx_2\; dt\right.$$
  \begin{eqnarray}\label{l1c}
  \ +\left.\int_0^T\int_{\mathcal D}\int^{L}_\alpha\int_{\alpha}^{x_1}  | {F(t,\xi,x_2)} | ^2 e^{-2s\eta(t,\xi,x_2)}r(x_1,\xi)d\xi\; dx_1\; dx_2\; dt\right].
  \end{eqnarray}
 For the first term  on the right hand side of (\ref{l1c}), formula (\ref{l1b}) implies
  \begin{eqnarray*}
  &\int_0^T\int_{\mathcal D}\int_{-L}^\alpha\int^{\alpha}_{x_1}  | {F(t,\xi,x_2)} | ^2 e^{-2s\eta(t,\xi,x_2)}r(x_1,\xi)d\xi\; dx_1\; dx_2\; dt\\
  &\leq \int_0^T\int_{\mathcal D}\int_{-L}^\alpha\int^{\alpha}_{x_1}  | {F(t,\xi,x_2)} | ^2 e^{-2s\eta(t,\xi,x_2)}d\xi\; dx_1\; dx_2\; dt
  \end{eqnarray*}
  and for   the second term on the right hand side of (\ref{l1c}) we obtain
  \begin{eqnarray*}
  &\int_0^T\int_{\mathcal D}\int^{L}_\alpha\int_{\alpha}^{x_1}  | {F(t,\xi,x_2)} | ^2 e^{-2s\eta(t,\xi,x_2)}r(x_1,\xi)d\xi\; dx_1\; dx_2\; dt\\
 &\leq \int_0^T\int_{\mathcal D}\int^{L}_\alpha\int_{\alpha}^{x_1}  | {F(t,\xi,x_2)} | ^2 e^{-2s\eta(t,\xi,x_2)}d\xi\; dx_1\; dx_2\; dt.
 \end{eqnarray*}
 We deduce easily Lemma \ref{l1} from these estimates.
 \begin{flushright}
\rule{.05in}{.05in}
\end{flushright}

  \noindent
  
  \noindent
  Now let use return to the stability result.
  %
  \subsection{Stability result}
  %
  \noindent
In this subsection we consider $u$ $\widetilde u$, $v$, $w$ and $z$ introduced in the previous subsection. Our goal is to use the Carleman estimate (\ref{ca2}). For this purpose, we will exploit the fact that in the bounded wave guide $\Omega$ a derivation with respect to $x_1$ does not alter the Dirichlet condition on $ \Sigma_1$ and that the Neumann condition on $\Sigma_2$ becomes a Dirichlet condition. 
The main result of this subsection is  the following stability estimate.
\begin{theorem}\label{t4}
Let Assumptions  \ref{as1}  be fulfilled. Let  $r>0$ be such that $$r\geq\textrm{max}\left(\ | {q}\ | _{L^2\left((0,T)\times\mathcal{D}\right)},\ | {\widetilde q}\ | _{L^2\left((0,T)\times\mathcal{D}\right)}\right).$$ Then, for any $0<\varepsilon<\frac{T}{2}$, there exists a constant $C_\varepsilon>0$ depending of $\varepsilon$, $b$, $u_0$, $k^\pm$ and $r$ such that
%
$$\| {q-\widetilde q}\| ^2_{L^2\left((\varepsilon,T-\varepsilon)\times\mathcal{D}\right)} \leq C_\varepsilon\left[\| {\partial_\nu\partial_{x_1}\widetilde u-\partial_\nu\partial_{x_1}u}\| ^2_{L^2((0,T)\times\Gamma^+)} \right.$$
$$\left.+ \| {\widetilde u(\cdot,\alpha,\cdot)-u(\cdot,\alpha,\cdot)}\| ^2_{H^1_t\left(0,T,\  H^2_{x_2}\left(\mathcal{D}\right)\right)}\right]$$

%
\end{theorem}
{\bf{Proof}}\\
According to \cite{LSU} and the maximum principle, Assumptions \ref{as1} implies that $u,\widetilde u, \partial_{x_2}u, \partial_{x_2}\widetilde u\in\mathcal C^1(Q)$ and  $f\widetilde u\geq c_1>0$. Thus, $w,\ \partial_{x_2}w\in\mathcal C^1(Q)$ and one can write
\begin{equation}\label{taa2}w(t,x)=\int_\alpha^{x_1}z(t,x',x_2)\ d x'+w(t,\alpha,x_2),\quad (t,x_1,x_2)\in Q,\end{equation}
\begin{equation}\label{taa3}\partial_{x_2}w(t,x)=\int_\alpha^{x_1}\partial_{x_2}z(t,x',x_2)\ d x'+\partial_{x_2}w(t,\alpha,x_2),\quad (t,x_1,x_2)\in Q.\end{equation}
Let us consider the source term $B_2\partial_{x_2}w+bw$ of (\ref{syst1-z}). 
Using representations (\ref{taa2}) and (\ref{taa3}), we get
\begin{eqnarray*}
\int_Qe^{-2s\eta}(B_2\partial_{x_2}w+bw)^2\ d x \ d t \leq4\int_Qe^{-2s\eta}\left(B_2\int_\alpha^{x_1}\partial_{x_2}z(t,x',x_2)\ d x'\right)^2\ d xÊ \ d t\\
+4\int_Qe^{-2s\eta}\left(B_2\partial_{x_2}w(t,\alpha,x_2)\right)^2\ d x\ d t
+4\int_Qe^{-2s\eta}\left(b\int_\alpha^{x_1}z(t,x',x_2)\ d x'\right)^2\ d x\  d t\\
+4\int_Qe^{-2s\eta}\left(bw(t,\alpha,x_2)\right)^2\ d x\ \ d t.
\end{eqnarray*}
Then, applying Lemma \ref{l1}, we obtain 
$$
\int_Qe^{-2s\eta}(B_2\partial_{x_2}w+bw)^2\ d x \ d t\leq C\left(\int_Qe^{-2s\eta}\left( | {\partial_{x_2}(z)} | ^2+ | z | ^2\right)\ dx \ dt\right)$$
\begin{eqnarray}\label{taab}
\ +C\left(\int_Qe^{-2s\eta}\left( | w(t,\alpha,x_2) | ^2+ | {\partial_{x_2}w(t,\alpha,x_2)} | ^2\right)\ dx\ dt\right).
\end{eqnarray}
Note that Assumption \ref{psi} implies
$$\left(\int_Qe^{-2s\eta}\left( | {w(t,\alpha,x_2)} | ^2+ | {\partial_{x_2}w(t,\alpha,x_2)} | ^2\right)\ dx\ dt\right)$$
$$\leq C\int_0^T\int_{\mathcal{D}} (| {w(t,\alpha,x_2)} | ^2+ | {\partial_{x_2}w(t,\alpha,x_2)} | ^2)\ dx_2\  dt.$$
Combining this estimate with (\ref{taab}), we obtain

$$\int_Qe^{-2s\eta}(B_2\partial_{x_2}w+bw)^2\ dx \ dt\leq C\left(\int_Qe^{-2s\eta}\left( | {\partial_{x_2}(z)} | ^2+ | {z} | ^2\right)\ dx \ dt\right.$$
\begin{eqnarray}\label{taad}\left.+\int_0^T\int_{\mathcal{D}}\left( | {w(t,\alpha,x_2)} | ^2+ | {\partial_{x_2}w(t,\alpha,x_2)} | ^2\right)\ dx_2\ dt\right).
\end{eqnarray}
An application of the Carleman estimate (\ref{ca2}) to $z$ yields
  \begin{eqnarray*}
  \int_Qe^{-2s\eta}\left[(s g)^{-1}(\Delta z)^2+(s g)^{-1}(\partial_t z)^2+s g\ | {\nabla z} | ^2+(sg)^3 z^2\right]\ dx\ dt\\ \nonumber
\leq C\left(\int_Qe^{-2s\eta}(B_2\partial_{x_2}w+bw)^2\ dx\  dt+\int_{(0,T)\times\Gamma^+}e^{-2s \eta}s g(\partial_\nu z)^2\ d\sigma \ d t\right).
  \end{eqnarray*}
Combining this estimate with (\ref{taad}), we obtain
  \begin{eqnarray*}\int_Qe^{-2s\eta}\left[(s g)^{-1}(\Delta z)^2+(s g)^{-1}(\partial_t z)^2+s g | {\nabla z} | ^2+(sg)^3 z^2\right]\ d x\ d t\\
 \leq C\left(\int_Qe^{-2s\eta}\left( | {\partial_{x_2}(z)} | ^2+ | {z} | ^2\right)\ d x \ d t
+\int_0^T\int_{\mathcal{D}}\left( | {w(t,\alpha,x_2)} | ^2+ | {\partial_{x_2}w(t,\alpha,x_2)} | ^2\right)\ d x_2\ d t\right.\\ 
\left.+\int_{(0,T)\times\Gamma^+}e^{-2s \eta}s g(\partial_\nu z)^2\ d\sigma \ d t\right).
  \end{eqnarray*}
Then,  for  $s$ sufficiently large, we get
\begin{eqnarray}\label{taa5}
\int_{Q_\varepsilon} \left( | {\partial_tz} | ^2+|\Delta z|^2+|\nabla z|^2+| z|^2\right)\ d
x  \ d t\\ 
\leq C_\varepsilon\left(\int_{0}^T  
\int_{\Gamma^+}  |\partial_{\nu} z|^2\ \ d \sigma \ d t
+\int_0^T\int_{\mathcal D}\left( | {w(t,\alpha,x_2)} | ^2+ | {\partial_{x_2}w(t,\alpha,x_2)} | ^2\right)\ d x_2\ d t\right).\nonumber
\end{eqnarray}
 with $Q_\varepsilon=(\varepsilon,T-\varepsilon)\times (-L,L)\times\mathcal{D}$. 
 Now, note that \[\widetilde{q}(t,x_2)-q(t,x_2)=\partial_t w -\Delta w+\mathbb{A}\cdot \nabla w +aw=Pw.\]
 Thus, applying (\ref{taa2}), one get
\[\widetilde{q}(t,x_2)-q(t,x_2)=P\int_\alpha^{x_1}z(t,x',x_2)\ d x'+Pw(t,\alpha,x_2).\]
Then, this representation  and estimate (\ref{taa5}) imply
 \begin{eqnarray*}
 \int_\varepsilon^{T-\varepsilon}\int_{\mathcal{D}}\left(\widetilde{q}(t,x_2)-q(t,x_2)\right)^2\ d x_2\ d t \\ \nonumber \leq C\left(\int_{Q_\varepsilon} \left( | {\partial_tz} | ^2+|\Delta z|^2+|\nabla z|^2+| z|^2\right)\ dx  \ d t  
 +\| {w(\cdot,\alpha,\cdot)}\ | ^2_{H_t^1\left(0,T, H_{x_2}^2(\mathcal{D})\right)}\right)\\ \nonumber 
 \leq C_\varepsilon\left(\int_{0}^T  
\int_{\Gamma^+}  |\partial_{\nu} z|^2\ \ d \sigma \ d t
+ \ | {w(\cdot,\alpha,\cdot)}\ | ^2_{H_t^1\left(0,T, H_{x_2}^2(\mathcal{D})\right)}\right).
\end{eqnarray*}
This completes the proof.
 \begin{flushright}
\rule{.05in}{.05in}
\end{flushright}

\section{Stability for an unbounded waveguide}
%
In this section we consider problem (\ref{systu}). We will establish a stability
result and deduce a uniqueness  for the
coefficient $q$. We give the result for the open waveguide with Dirichlet boundary   conditions. The Carleman estimate for unbounded domain will be the key
ingredient in the proof of such a  stability estimate.\\
\noindent
 From now on, we set \[\Omega=\mathbb{R}\times \mathcal{D},\  Q=(0,T)\times\Omega,\  \Sigma=(0,T)\times\partial \Omega.\] 
\subsection{Global Carleman Estimate}
%
%
\setcounter{equation}{0}
Let $f(x_1)$ be a bounded positive function in $\mathcal{C}^2(\mathbb{R})$ such that
$f(x_1) \ge c_{min}>0$, $f$ and all its derivatives up to order two are 
bounded by a positive constant $\tilde{C_0}.$\\
\noindent
Let $u=u(t,x)$ be a function equals to zero on $(0,T)\times\partial \Omega$ and solution
of the Heat equation 
$$\partial_t u -\Delta u=F.$$
We prove here a global Carleman-type estimate for $u$
with a single observation 
acting on a part $\Gamma^+_1$ of the boundary $\Gamma$ of $\mathcal{D}$
in the right-hand side of the estimate.  
Let $\psi_2$ be the function defined in subsection 2.1 and let  ${\psi}$ be a 
${\cal C}^4( \mathbb{R}^d)$ function defined by
\[\psi(x_1,x_2)=e^{x_1}\psi_2(x_2).\]
Then, the function $\psi$ satisfies the conditions:
\begin{assumption}\label{funct-psi}
\begin{itemize}
\item $\psi (x) >0$ in $\overline{\Omega}$,
\item  There exists $C_0>0$ such that $|\nabla {\psi}| \geq C_0>0 \;\;\mbox{ in } \;\;\Omega$,
\item ${\partial}_{\nu} {{\psi}}\leq 0\;\;\mbox{on}\;\;(0,T)\times \mathbb{R} \times\Gamma_1^-$,
\item  $\inf_{x\in\overline{\Omega}}\partial_{x_1}\psi(x)>0,$
\item $\displaystyle \lim_{x_1\to\pm\infty}\inf_{x_2\in\overline{D}} \left\vert {\frac{\psi(x_1,x_2)}{x_1}} \right\vert =+\infty.$
\end{itemize}
\end{assumption}
Now let us introduce the function 
$$
\varphi(t,x)=g(t)e^{\lambda \psi(x)},\ \ \lambda>0\ \mbox{ with }\  g(t)=\frac{1}{t(t-T)}.
$$
Let $H$ be the operator defined by
\begin{equation} \label{H}
Hu:=\partial_t u-\Delta u   \;\mbox{ in }\; {Q}=\Omega \times (0,T).
\end{equation}
We set $w=e^{-s \varphi}u$, $M w = e^{-s \varphi} H(e^{s \varphi} w)$ for $s>0$
 and we introduce the following operators
\begin{equation} \label{M1}
  M_1 w : =-\Delta w-s^2 |\nabla \varphi |^2w-s\partial_t \varphi w ,
\end{equation}
\begin{equation} \label{M2}
  M_2 w  : =\partial_t w +2s \nabla \varphi \cdot \nabla w 
 +s\Delta \varphi w .
\end{equation}
Then the following result holds.\\
\begin{theorem}
\label{th-Carl} 
Let $H$, $M_1$, $M_2$ be the operators defined respectively by
(\ref{H}), (\ref{M1}), (\ref{M2}). We assume that Assumptions \ref{funct-psi}
are satisfied. 
Then there exist $\lambda_0> 0$, $s_0>0$ and a positive
constant $C=C(\Omega, \Gamma,T)$ such that, for any $\lambda \ge
\lambda_0$ and any $s \ge s_0 $, the next inequality holds:
\begin{eqnarray}
\label{Carl}
   s^{3} \lambda^{4} \int_{Q} e^{-2s \varphi} \varphi^{3} |u|^2\ d x  \ d t  
+s \lambda    \int_{Q} e^{-2s \varphi} \varphi |\nabla u|^2\ d
x \ d t
  +  \|M_1(e^{-s \varphi}u)\|^2_{L^2({Q})} \nonumber\\
+ 
\|M_2(e^{-s \varphi}u)\|^2_{L^2({Q})}
  \leq C \left[
   s \lambda \int_{0}^T  \int_{\mathbb{R}_{x_1}}
\int_{\Gamma^+_1} e^{-2s \varphi} \varphi  |\partial_{\nu} u |^2\ \partial_{\nu} \psi \ d \sigma \ d t\right.\\
 \left. +  \int_{Q} e^{-2s \varphi}\ | H u |^2\ d x \ d t\right], \nonumber
\end{eqnarray}
for all $u$ satisfying $Hu\in L^2(\Omega \times (0,T)),$ $u\in
L^2(0,T;H^1_0(\Omega)),$ $\partial_{\nu} u\in L^2(0,T;L^2(\Gamma)).$
\end{theorem}
Note that this theorem has already been proved in \cite{FI}, \cite{I} and \cite{T}. 
\noindent
In the inequality (\ref{Carl}) , we can also have an estimate of $\partial_t u$ and $\Delta u$ (see \cite{FI}).
%
\subsection{Inverse Problem}
%
In this subsection, we establish a stability
result and deduce a uniqueness result for the
coefficient $q$. \\
The Carleman estimate (\ref{Carl}) will be the key
ingredient in the proof of such a  stability estimate.\\
Let $u$ be solution of
$$\left \{ \begin{array}{lll}
   \partial_t u -\Delta u+q(t,x_2)f(x_1)u=0 & \mbox{in} &Q,\\
   u(t,x)=b(t,x)  & \mbox{on} & \Sigma,\\
   u(0,x)=u_0(x) & \mbox{in} & \Omega,
\end{array}\right.$$
and $\widetilde{u}$ be solution of
$$\left \{ \begin{array}{lll}
  \partial_t \widetilde{u} -\Delta \widetilde{u}+\widetilde{q}(t,x_2)f(x_1)\widetilde{u}  =0 & \mbox{in} & Q,\\
   \widetilde{u}(t,x)=b(t,x)  & \mbox{on} &\Sigma,\\
   \widetilde{u}(0,x)=u_0(x) & \mbox{in} & \Omega,
\end{array}\right.$$
%
\begin{assumption} \label{g} Here we assume that: 
\begin{itemize}
\item $q(t,x_2)f(x_1), \widetilde{q}(t,x_2)f(x_1)\in\mathcal C^{1+\alpha,1+\frac{\alpha}{2}}(\overline{Q})\cap L^\infty (Q),$
\item $b\in\mathcal C^{3+\alpha,1+\frac{\alpha}{2}}(\overline{\Sigma})\cap L^2((0,T);H^{\frac{5}{2}}(\partial\Omega))\cap H^{\frac{3}{4}}((0,T);H^{1}(\partial\Omega)),$
\item $u_0\in \mathcal C^{3,\alpha}(\overline{\Omega})\cap H^3(\Omega),$
\item $\partial_tb(0,x)-\Delta_xu_0(x)+q(0,x_2)f(x_1)u_0(x)=0,\quad x=(x_1,x_2)\in\partial\Omega,$
\item $\partial_tb(0,x)-\Delta_xu_0(x)+\widetilde{q}(0,x_2)f(x_1)u_0(x)=0,\quad x=(x_1,x_2)\in\partial\Omega,$
\item There exists $r>0$ such that $b\geq r$ and $u_0\geq r$.
\end{itemize}
\end{assumption}

If we set $v=u- \widetilde{u}$, then v satisfies 
%
$$\left \{ \begin{array}{lll}
 \partial_t v -\Delta v+q(t,x_2)f(x_1)v=(\widetilde{q}(t,x_2)-q(t,x_2))f(x_1)\widetilde{u}
& \mbox{in} & Q,\\
   v(t,x)=0  & \mbox{on} & \Sigma,\\
   v(0,x)=0 & \mbox{in} & \Omega.
\end{array}\right.$$
%
Then with the change of function $\displaystyle{w=\frac{v}{f\widetilde{u}}}$, $w$ is solution of the following system
$$\left \{ \begin{array}{lll}
   \partial_t w -\Delta w+\mathbb{A}\cdot \nabla w +aw=\widetilde{q}(t,x_2)-q(t,x_2)
 & \mbox{in} &Q,\\
   w(t,x)=0  & \mbox{on} &\Sigma,\\
   w(0,x)=0 & \mbox{in} & \Omega.
\end{array}\right.$$
%
where
$$\mathbb{A}=\frac{-2}{f\widetilde{u}}\nabla(f\widetilde{u}) \ \ \mbox{and}\ \ 
a=\frac{\partial_t(f\widetilde{u})-\Delta(f\widetilde{u})+V}{f\widetilde{u}}.$$
We consider the $x_1$-derivative of the previous system and we set $z:=\partial_{x_1} w$, then
$z$ is solution of
\begin{equation} \label{syst-z}
\left \{ \begin{array}{lll}
   \partial_t z -\Delta z+\mathbb{A}\cdot \nabla z +az+B_1z+B_2\partial_{x_2}w +bw=0
    & \mbox{in} & Q,\\
   z(t,x)=0  & \mbox{on} & \Sigma,\\
   z(0,x)=d(x) & \mbox{in} & \Omega
\end{array}\right.
\end{equation}
with 
$$B_1:=-2\partial_{x_1}\left(\frac{\partial_{x_1}(f\widetilde{u})}{f\widetilde{u}}\right),\ \ B_2:=-2\partial_{x_1}\left(\frac{\partial_{x_2}(f\widetilde{u})}{f\widetilde{u}}\right)    \ \ \mbox{and}\ \ b=\partial_{x_1} a.$$
\begin{assumption}\label{h} 
$\mathbb{A}$, $B_1$, $B_2$ and $a$, $b$ are bounded.
\end{assumption}
We can apply the Carleman estimate (\ref{Carl}) for $z$ and we obtain:
\begin{eqnarray} \label{Iz}
I(z) \leq C\left[
   s \lambda \int_{0}^T   \int_{\mathbb{R}_{x_1}}
\int_{\Gamma^+_1} e^{-2s\varphi} \varphi  |\partial_{\nu}  z|^2\ \partial_{\nu} \psi \ d \sigma \ d t \right.\\ \nonumber
 \left.   \int_{Q} e^{-2s\varphi}\ (| \nabla z |^2+ | z |^2)\ d x \ d t+
 \int_{Q} e^{-2s\varphi}\ (| \nabla w |^2+ | w |^2)\ d x \ d t\right].
  \end{eqnarray}
  The second integral of the right hand side of the previous estimate is "absorbed"
  by the left hand side, for $s$ sufficiently large. For the last integral, we need the following lemma which is an adaptation of a lemma proved in \cite{K} and \cite{KT}:
  \begin{lemma} \label{lemKli}
  Let $F$ be a function in $L^2(Q)\cap \mathcal C (\overline{Q})$. Then we have the following estimate:
 $$ \int_Q \left |  \int_{\alpha}^{x_1} F(t,\xi,x_2) d\xi \right | ^2 e^{-2s\varphi} dx_1\; dx_2\; dt \leq \frac{C}{s^2}
   \int_Q \left |  F(t,x)  \right | ^2 e^{-2s\varphi} dx\; dt .$$
  \end{lemma}
 { \bf{Proof}}\\
 We recall here the proof of this lemma.
  \begin{eqnarray*}
  I & := &  \int_Q \left |  \int_{\alpha}^{x_1} F(t,\xi,x_2) d\xi \right | ^2 e^{-2s\varphi} dx_1\; dx_2\; dt\\
  & = &  \int_0^T \int_{\mathcal{D}} \int_{-\infty}^{\alpha} \left |  \int_{\alpha}^{x_1} F(t,\xi,x_2) d\xi \right | ^2 e^{-2s\varphi} dx_1\; dx_2\; dt \\
  &&+
  \int_0^T \int_{\mathcal{D}} \int_{\alpha}^ {+\infty}\left |  \int_{\alpha}^{x_1} F(t,\xi,x_2) d\xi \right | ^2 e^{-2s\varphi} dx_1\; dx_2\; dt\\
  &:=& I_1+I_2.
  \end{eqnarray*}
  Note that
  $$e^{-2s\varphi}=\frac{2s \partial_{x_1} \varphi}{2s \partial_{x_1} \varphi}e^{-2s\varphi}=
  \frac{1}{2s \partial_{x_1} \varphi}\partial_{x_1}(-e^{_2s\varphi})$$
  and according to Assumption \ref{funct-psi} there exists a positive constant $\kappa$ such that \[\partial_{x_1}\varphi \geq \kappa>0.\] So, we have
  $$e^{-2s\varphi}=
  \frac{1}{2s \partial_{x_1} \varphi}\partial_{x_1}(-e^{-2s\varphi})
  \leq \frac{1}{2s \kappa}\partial_{x_1}(-e^{-2s\varphi}).$$
  We first give an estimate for $I_1$:
   \begin{eqnarray*}
  I_1 & = &  \int_0^T \int_{\mathcal{D}} \int_{-\infty}^{\alpha} \left |  \int_{\alpha}^{x_1} F(t,\xi,x_2) d\xi \right | ^2 e^{-2s\varphi} dx_1\; dx_2\; dt\\
  &\leq &   \frac{1}{2s \kappa}\int_0^T \int_{\mathcal{D}} \int_{-\infty}^{\alpha} \left |  \int_{\alpha}^{x_1} F(t,\xi,x_2) d\xi \right | ^2 \partial_{x_1}(-e^{-2s\varphi}) dx_1\; dx_2\; dt.
  \end{eqnarray*}
  Using integration by parts and the Cauchy-Schwarz inequality, we obtain
  \begin{eqnarray*}
  I_1 & = &  \frac{1}{s \kappa}\int_Q e^{-s\varphi} F(t,x) \left( \int_{\alpha}^{x_1} F(t,\xi,x_2) d\xi \right )  e^{-s\varphi}dx_1\; dx_2\; dt\\
  &\leq & \frac{C}{2s \kappa}\left(\int_Q \left |  \int_{\alpha}^{x_1} F(t,\xi,x_2) d\xi \right | ^2 e^{-2s\varphi} dx_1\; dx_2\; dt\right)^{1/2} \left( \int_Q \left |  F(t,x)  \right | ^2 e^{-2s\varphi} dx\; dt\right)^{1/2}.
  \end{eqnarray*}
  That is
  $$I_1 \leq  \frac{C}{s \kappa} I^{1/2}  \left( \int_Q \left |  F(t,x)  \right | ^2 e^{-2s\varphi} dx\; dt\right)^{1/2}.$$
In the same way, we obtain 
  $$I_2 \leq  \frac{C}{s \kappa} I^{1/2}  \left( \int_Q \left |  F(t,x)  \right | ^2 e^{-2s\varphi} dx\; dt\right)^{1/2}$$
  and therefore 
  $$I \leq  \frac{2C}{s \kappa} I^{1/2}  \left( \int_Q \left |  F(t,x)  \right | ^2 e^{-2s\varphi} dx\; dt\right)^{1/2}.$$
 This complete the proof.
 \begin{flushright}
\rule{.05in}{.05in}
\end{flushright}
  \noindent
Now we come back to the inequality (\ref{Iz}) in order to estimate the last integral of the right hand side.\\
A direct application of Lemma \ref{lemKli} leads to:
 $$   \int_{Q} e^{-2s\varphi}\ (| \nabla w |^2+ | w |^2)\ d x  \ d t
 \leq \frac{C}{s^2}     \int_{Q} e^{-2s\varphi}\ (| \nabla z |^2+ | z |^2)\ d x  \ d t.$$
 So (\ref{Iz}) becomes
 \begin{eqnarray*} 
\hspace{-2cm}I(z) \leq C
   s \lambda \int_{0}^T   \int_{\mathbb{R}_{x_1}}
\int_{\Gamma^+_1} e^{-2s\varphi} \varphi  |\partial_{\nu}  z|^2\ \partial_{\nu} \psi \ d \sigma \ d t 
  +
 \frac{C}{s^2}   \int_{Q} e^{-2s\varphi}\ (| \nabla z |^2+ | z |^2)\ d x  \ d t.
  \end{eqnarray*}
  The last integral of the right hand side is "absorbed"
  by the left hand side, for $s$ sufficiently large. We thus obtain a Carleman estimate for $z$
  solution of (\ref{syst-z}):
$$I(z) \leq C
   s \lambda \int_{0}^T   \int_{\mathbb{R}_{x_1}}
\int_{\Gamma^+_1} e^{-2s\varphi} \varphi  |\partial_{\nu}  z|^2\ \partial_{\nu} \psi \ d \sigma \ d t $$
\subsection{Stability Estimate}

\noindent
In this subsection we consider $u$ $\widetilde u$, $v$, $w$ and $z$ introduced in the previous subsection. We will exploit the fact that in the wave guide $\Omega$ derivations with respect to $x_1$ do not alter the Dirichlet condition. 
The main result of this subsection is  the following stability estimate.
\begin{theorem}\label{t3}
Let Assumptions   \ref{g} and \ref{h} be fulfilled. Let $\alpha\in\mathbb{R}$ and $r>0$ be such that $r\geq\textrm{max}\left(\|{q}\|_{L^2\left((0,T)\times\mathcal{D}\right)},\|{\widetilde q}\|_{L^2\left((0,T)\times\mathcal{D}\right)}\right)$. Then, for any $0<\varepsilon<\frac{T}{2}$, there exists a constant $C_\varepsilon>0$ depending of $\varepsilon$, $b$, $u_0$, $k^\pm$ and $r$ such that
\begin{eqnarray} \label{ta1}
\|{q-\widetilde q}\|^2_{L^2\left((\varepsilon,T-\varepsilon)\times\mathcal{D}\right)} &\leq& C_\varepsilon\left[\|{\partial_\nu\partial_{x_1}\widetilde u-\partial_\nu\partial_{x_1}u}\|^2_{L^2((0,T)\times\mathbb{R}_{x_1}\times\Gamma^+_1)} \right.\\ \nonumber
&&\left.+ \|{\widetilde u(\cdot,\alpha,\cdot)-u(\cdot,\alpha,\cdot)}\|^2_{H^1_t\left(0,T,\  H^2_{x_2}\left(\mathcal{D}\right)\right)}\right]
\end{eqnarray}

\end{theorem}
{\bf{Proof}}\\
According to \cite{LSU} and the maximum principle, Assumptions \ref{g} implies that $u,\widetilde u, \partial_{x_2}u, \partial_{x_2}\widetilde u\in\mathcal C^1(Q)$ and  $f\widetilde u\geq c_1>0$. Thus, $w,\ \partial_{x_2}w\in\mathcal C^1(Q)$ and one can write
$$\label{ta2}w(t,x)=\int_\alpha^{x_1}z(t,x',x_2)\ dx'+w(t,\alpha,x_2),\quad (t,x_1,x_2)\in Q,$$
$$\label{ta3}\partial_{x_2}w(t,x)=\int_\alpha^{x_1}\partial_{x_2}z(t,x',x_2)\ d x'+\partial_{x_2}w(t,\alpha,x_2),\quad (t,x_1,x_2)\in Q.$$
Let us consider the source term $B_2\partial_{x_2}w+bw$ of (\ref{syst-z}). 
Combining  Lemma \ref{lemKli} with some arguments used in Theorem \ref{t4}, we obtain 
\begin{eqnarray}\label{tab}
\hspace{-2,5cm}\int_Qe^{-2s\varphi}(B_2\partial_{x_2}w+bw)^2\ d x \ d t\leq& \frac{C}{s^2}\left(\int_Qe^{-2s\varphi}\left(|{\partial_{x_2}(z)}|^2+|{z}|^2\right)\ d x \ d t\right.	\\
\ &\left.\ \ +\int_Qe^{-2s\varphi}\left(|{w(t,\alpha,x_2)}|^2+|{\partial_{x_2}w(t,\alpha,x_2)}|^2\right)\ d x\ d t\right).\nonumber
\end{eqnarray}
Since
 $$\sup_{(t,x_2)\in(0,T)\times\mathcal{D}}\int_\mathbb{R} e^{-s\varphi(t,x_1,x_2)} dx_1<+\infty$$ 
we deduce that
  \begin{eqnarray*}\left(\int_Qe^{-s\varphi}\left(|{w(t,\alpha,x_2)}|^2+|{\partial_{x_2}w(t,\alpha,x_2)}|^2\right)\ d x\ d t\right)\nonumber\\
  \leq C\left(\int_0^T\int_{\mathcal{D}}\left(|{w(t,\alpha,x_2)}|^2+|{\partial_{x_2}w(t,\alpha,x_2)}|^2\right)\ d x_2\ d t\right).\end{eqnarray*}
Combining this estimate with (\ref{tab}), we obtain
\begin{eqnarray}\label{tad}
\hspace{-2cm}\int_Qe^{-2s\varphi}(B_2\partial_{x_2}w+bw)^2\ d x \ d t\leq& \frac{C}{s^2}\left(\int_Qe^{-2s\varphi}\left(|{\partial_{x_2}(z)}|^2+|{z}|^2\right)\ d x \ d t\right.	\\
\ &\left.+\int_0^T\int_{\mathcal{D}}\left(|{w(t,\alpha,x_2)}|^2+|{\partial_{x_2}w(t,\alpha,x_2)}|^2\right)\ d x_2\ d t\right).\nonumber
\end{eqnarray}
An application of the Carleman estimate (\ref{Carl}) to $z$ yields
$$s^{3} \lambda^{4} \int_{Q} e^{-2s\varphi} \varphi^{3} |z|^2\ d x  \ d t  
+s \lambda \int_{Q} e^{-2s\varphi} \varphi |\nabla z|^2\ d
x  \ d t
  +  \|M_1(e^{-s \varphi}z)\|^2_{L^2({Q})}+ 
\|M_2(e^{-s \varphi}z)\|^2_{L^2({Q})}$$
$$
  \leq C \left[
   s \lambda \int_{0}^T  \int_{\mathbb{R}_{x_1}}
\int_{\Gamma^+_1} e^{-2s\varphi} \varphi  |\partial_{\nu} z|^2\ \partial_{\nu} \psi \ d \sigma \ d t
  +  \int_{Q} e^{-2s\varphi}\ | B_2\partial_{x_2}w+bw |^2\ d x  \ d t\right].
$$
Combining this estimate with (\ref{tad}), we find
$$s^{3} \lambda^{4}  \int_{Q} e^{-2s\varphi} \varphi^{3} |z|^2\ d x  \ d t  
+s \lambda   \int_{Q} e^{-2s\varphi} \varphi |\nabla z|^2\ d
x  \ d t 
+ s^{-1}  \int_{Q} e^{-2s\varphi} \varphi^{-1} \left(|{\partial_tz}|^2+|\Delta z|^2\right)\ d
x  \ d t $$
$$\leq C 
   s \lambda \int_{0}^T   \int_{\mathbb{R}_{x_1}}
\int_{\Gamma^+_1} e^{-2s\varphi} \varphi  |\partial_{\nu} z|^2\ \partial_{\nu} \psi \ d \sigma \ d t
  +\frac{C}{s^2}\left(\int_Qe^{-2s\varphi}\left(|{\partial_{x_2}(z)}|^2+|{z}|^2\right)\ d x \ d t\right)$$
$$  +C\left(\int_0^T\int_{\mathcal D}\left(|{w(t,\alpha,x_2)}|^2+|{\partial_{x_2}w(t,\alpha,x_2)}|^2\right)\ d x_2\ d t\right).\nonumber
$$
Then, using Assumption \ref{funct-psi}, for $s$ and $\lambda$ sufficiently large we get
$$
\int_{Q_\varepsilon} \left(|{\partial_tz}|^2+|\Delta z|^2+|\nabla z|^2+| z|^2\right)\ d
x  \ d t 
\leq C\left(\int_{0}^T  
\int_{\Gamma^+_1} \varphi  |\partial_{\nu} z|^2\ \partial_{\nu} \psi \ d \sigma \ d t\right)$$
\begin{eqnarray}\label{ta5}
+C_\varepsilon\left(\int_0^T\int_{\mathcal D}\left(|{w(t,\alpha,x_2)}|^2+|{\partial_{x_2}w(t,\alpha,x_2)}|^2\right)\ d x_2\ d t\right).
\end{eqnarray}
 with $Q_\varepsilon=(\varepsilon,T-\varepsilon)\times (-R,R)\times\mathcal{D}$. 
Combining this estimate with some arguments used in Theorem \ref{t4} , we deduce easily (\ref{ta1}) and the proof of Theorem \ref{t3} is complete.
 \begin{flushright}
\rule{.05in}{.05in}
\end{flushright}
%

\section*{References}
%
   
%

\begin{thebibliography}{99}
%
\bibitem[AT]{AT}{\sc P. Albano and D. Tataru}, {\em Carleman estimates and boundary observability 
for a coupled parabolic-hyperbolic system}. Electronic Journal of Differential Equations, vol 2000, 22, 1--15, (2000).
%
\bibitem[BK]{BK}{\sc A.L. Bukhgeim and M.V. Klibanov}, {\em Uniqueness in the large of a class of multidimensional
inverse problems}, Soviet Math. Dokl., 17, 244-247 (1981).
%
\bibitem[CE]{CE}{\sc J. R. Cannon and S. P. Esteva}, {\em An inverse problem for the heat equation}, Inverse Problems, {\bf2} (1986), 395-403.
%
\bibitem[Ch09]{Ch09}{\sc M. Choulli}, {\em Une introduction aux probl\`emes inverses elliptiques et paraboliques}, Math\'ematiques et Applications, Vol. 65, Springer-Verlag, Berlin, 2009.
%
\bibitem[CIK]{CIK}{\sc D. Chae, O. Yu Imanuvilov and S. M; Kim}, {\em Exact controllability for semi-linear parabolic equations
with Neumann boundary conditions}, J. Dynam. Control Systems 2 (4), 449 -- 483 (1996).
%
\bibitem[CY06]{CY06}{\sc M. Choulli and  M. Yamamoto},
{\em Some stability estimates in determining sources and coefficients},
J. Inv. Ill-Posed Problems, {\bf14} (4) (2006), 355-373.
%
\bibitem[CY11]{CY11}{\sc M. Choulli and  M. Yamamoto}, {\em Global existence and stability for an inverse coefficient problem for a semilinear parabolic equation}, Arch. Math (Basel),  {\bf97} (6) (2011), 587-597.
%
\bibitem[E07]{E07}{\sc G. Eskin}, {\em Inverse Hyperbolic Problems with Time-Dependent Coefficients}, Commun. PDE, {\bf32}  (11) (2007), 1737-1758.
%
\bibitem[E08]{E08}{\sc G. Eskin}, {\em Inverse problems for the Schr\"odinger equations with time-dependent electromagnetic potentials and the Aharonov-Bohm effect}, J. Math. Phys., {\bf49} (2) (2008), 1-18.
%
\bibitem[FI]{FI}{\sc A.V. Fursikov and O.Y. Imanuvilov}, {\em Controllability of evolution equations}, Lecture Notes Series,
34, Seoul National University, (1996).
%
\bibitem[I]{I}{\sc O. Yu. Immanuvilov}, {\em Controllability of parabolic equationss}, Sbornik Mathematics,
186, no 6, 879--900 (1993).
%
\bibitem[K]{K}{\sc M.V. Klibanov}, {\em Global uniqueness of 
a multidimensionnal inverse problem for a non linear parabolic 
equation by a Carleman estimate}, Inverse Problems, 20, 1003-1032, 
(2004).
%
\bibitem[KT]{KT}{\sc M.V. Klibanov and A. Timonov}, {\em Carleman estimates
for coefficient inverse problems and numerical applications}, Inverse and Ill-posed series,VSP, Utrecht (2004).
%
\bibitem[LSU]{LSU}{\sc O. A. Ladyzhenskaja, V. A. Solonnikov and N. N. Ural'tzeva}, {\em Linear and
quasilinear equations of parabolic type}, Nauka, Moscow, 1967 in Russian ;
English translation : American Math. Soc., Providence, RI, 1968.
%
 
%

%
\bibitem[T]{T}{\sc D. Tataru}, {\em Carleman estimates, Unique Continuation and Controllability for anisotropic
PDE's}, Contemporary Mathematics, vol 209, 267--279 (1997).
%
\bibitem[Y]{Y}{\sc M. Yamamoto}, {\em Carleman estimates for parabolic equations and applications},
Inverse Problems, 25 (2009).
%
%
\end{thebibliography}
\end{document}